\documentclass{amsart}

  \usepackage{amsmath}
  \usepackage{amssymb}
  \usepackage{amsthm}
  \usepackage{mathrsfs}
  \usepackage{bbm}
  \usepackage{pifont}
  \usepackage{eufrak}
  \usepackage{graphicx}

\def\ie{{\it i.e.\/} }
\def\TH{\mathop{T_H}}
\def\tn{\mathop{\tau_n}}
\def\tnp{\mathop{\tau_{n+1}}}
\def\NN{{\ifmmode{\mathbbm{N}}\else{$\mathbbm{N}$}\fi}}
\def\ZZ{{\ifmmode{\mathbbm{Z}}\else{$\mathbbm{Z}$}\fi}}
\def\A{{\ifmmode{\mathscr{A}}\else{$\mathscr{A}$}\fi}}
\def\H{{\ifmmode{\mathscr{H}}\else{$\mathscr{H}$}\fi}}
\def\B{{\ifmmode{\mathscr{B}}\else{$\mathscr{B}$}\fi}}
\def\build#1_#2^#3{\mathrel{\mathop{\kern 0pt#1}\limits_{#2}^{#3}}}
\def\converge#1#2#3{\build\hbox to 15mm{\rightarrowfill}_{#1\rightarrow #2}
     ^{\hbox{\scriptsize #3}}}

\def\egdef{:=}
\def\Id{\mathop{\mbox{Id}}}
  \newtheorem{theorem}{Theorem}[section]
  \newtheorem{lemma}{Lemma}[section]
  \newtheorem{defi}{Definition}[section]
    \newtheorem{prop}{Proposition}[section]
\begin{document}
\bibliographystyle{amsplain}

\title{An extension which is relatively twofold mixing but not threefold mixing}
\author{Thierry \sc {de la Rue}}
\address{Laboratoire de Math\'ematiques Rapha\"el Salem\\
	UMR 6085 CNRS -- Universit\'e de Rouen\\
	Site Colbert\\
	F76821 Mont-Saint-Aignan Cedex}
\email{thierry.de-la-rue@univ-rouen.fr}

\begin{abstract}
We give an example of a dynamical system which is mixing relative to one of its factors, 
but for which relative mixing of order three does not hold.
\end{abstract}

\maketitle

\section{Factors, extensions and relative mixing}

\subsection{Factors, extensions and Rokhlin cocycle }

We are interested in dynamical systems $(X,\A,\mu,T)$, where $T$ is an ergodic automorphism
of the Lebesgue space $(X,\A,\mu)$. We will often designate such a system by simply the symbol $T$.
A \emph{factor} of $T$ is a sub-$\sigma$-algebra $\H$ of $\A$ such that $\H=T^{-1}\H$.

The canonical example of a system with factor is given by the \emph{skew product}, constructed from 
a dynamical system $(X_H,\A_H,\mu_H,\TH)$ (called the \emph{base} of the skew product) and a measurable
map $x\longmapsto S_x$ from $X_H$ to the group of automorphisms of some Lebesgue space $(Y,\B,\nu)$ 
(such a map is called a \emph{Rokhlin cocycle}). The transformation is defined on the product space
$(X_H\times Y, \A_H\otimes\B, \mu_H\otimes\nu)$ by
$$ \tilde T(x,y)\ =\ (\TH x, \mathop{S_x} y). $$
In this context, the sub-$\sigma$-algebra $\A_H\otimes \{Y,\emptyset\}$ is clearly a factor of $\tilde T$.

Since the work of Abramov and Rokhlin \cite{abram1}, this kind of construction is known to be the general 
model for a system with factor: If $\H$ is  a factor of $T$, then there exists an 
isomorphism $\varphi$ between $T$ and a 
skew product $\tilde T$ constructed as above, with $\varphi(\H)=\A_H\otimes \{Y,\emptyset\}$.
In such a situation, we say that $T$ is an \emph{extension} of $\TH$.

\subsection{Mixing relative to a factor}

To understand precisely the way a factor is embedded in the dynamical system, one is led to study
the behaviour of the system \emph{relative to the factor}; to this end, relative properties are defined
which are generalizations of absolute properties of dynamical systems.
For example, one can define weak-mixing relative to a factor (see e.g. \cite{furst3}), or
the property of being a K-system relative to a factor \cite{rahe2}.

We are interested in this work in the property of being mixing relative to a factor.

\begin{defi}
  Let $\H$ be a factor of the system $(X,\A,\mu,T)$. $T$ is said 
  \emph{$\H$-relatively mixing} if
  \begin{equation}
    \label{eq:defmr}
    \forall A,B\in\A,\quad \mu\left(A\cap T^{-k}B|\H\right)-\mu(A|\H)\mu(T^{-k}B|\H) 
    \converge{k}{+\infty}{proba}\ 0.
  \end{equation}
\end{defi}

As for the absolute property of mixing, it is possible to define mixing relative to a factor of
any order $n\ge 2$. The property described by~(\ref{eq:defmr}) corresponds to relative mixing 
of order 2 (twofold relative mixing); for relative mixing of order 3 (threefold relative mixing),
(\ref{eq:defmr}) should be replaced by

  \begin{multline}
    \label{eq:defmr3}
    \forall A,B,C\in\A,\\
     \mu\left(A\cap T^{-j}B\cap T^{-k}C|\H\right)-\mu(A|\H)\mu(T^{-j}B|\H)\mu(T^{-k}C|\H) \converge{j, k-j}{+\infty}{proba}\ 0.
  \end{multline}

Whether (absolute) twofold mixing implies threefold mixing is a well-known open problem in ergodic
theory. The main goal of this work is to show that as far as relative mixing is concerned, twofold
does not necessarily imply threefold.

\begin{theorem}
  \label{erm23}
  We can construct a dynamical system $(X,\A,\mu,T)$ with a factor $\H$ such that 
  $T$ is $\H$-relatively twofold mixing  but not $\H$-relatively threefold mixing.  
\end{theorem}

\section{An extension which is relatively twofold mixing  but not relatively threefold mixing.}

\subsection{The base}

\begin{figure}[htbp]
  \begin{picture}(0,0)%
\includegraphics{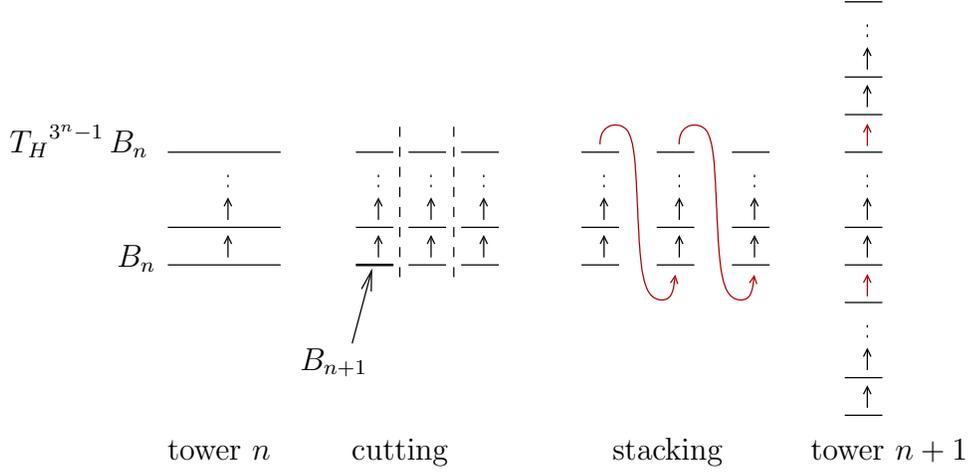}%
\end{picture}%
\setlength{\unitlength}{4144sp}%
\begingroup\makeatletter\ifx\SetFigFont\undefined%
\gdef\SetFigFont#1#2#3#4#5{%
  \reset@font\fontsize{#1}{#2pt}%
  \fontfamily{#3}\fontseries{#4}\fontshape{#5}%
  \selectfont}%
\fi\endgroup%
\begin{picture}(5232,2809)(-494,-1733)
\put(149,-524){\makebox(0,0)[lb]{\smash{\SetFigFont{12}{14.4}{\rmdefault}{\mddefault}{\updefault}$B_n$}}}
\put(450,-1679){\makebox(0,0)[lb]{\smash{\SetFigFont{12}{14.4}{\rmdefault}{\mddefault}{\updefault}tower $n$}}}
\put(1553,-1679){\makebox(0,0)[lb]{\smash{\SetFigFont{12}{14.4}{\rmdefault}{\mddefault}{\updefault}cutting}}}
\put(3113,-1679){\makebox(0,0)[lb]{\smash{\SetFigFont{12}{14.4}{\rmdefault}{\mddefault}{\updefault}stacking}}}
\put(4298,-1679){\makebox(0,0)[lb]{\smash{\SetFigFont{12}{14.4}{\rmdefault}{\mddefault}{\updefault}tower $n+1$}}}
\put(1243,-1140){\makebox(0,0)[lb]{\smash{\SetFigFont{12}{14.4}{\rmdefault}{\mddefault}{\updefault}$B_{n+1}$}}}
\put(-494,164){\makebox(0,0)[lb]{\smash{\SetFigFont{12}{14.4}{\rmdefault}{\mddefault}{\updefault}$\mathop{T_H}^{3^n-1}B_n$}}}
\end{picture}
  \caption{Construction of $\TH$ by cutting and stacking}
  \label{fig:base}
\end{figure}

The dynamical system announced in Theorem~\ref{erm23} is constructed as a skew product, 
whose base $(X_H,\A_H,\mu_H,\TH)$ is obtained as follows: Take $X_H\egdef [0,1[$ equipped with the
Lebesgue measure $\mu_H$ on the Borel $\sigma$-algebra $\A_H$. The transformation $\TH$ can be viewed
as a triadic version of the Von Neumann-Kakutani transformation; we describe now its construction 
by the \emph{cutting and stacking} method (see Figure~\ref{fig:base}).
 
We begin by splitting $X_H$ into three subintervals of length $1/3$; we set $B_1\egdef  [0,1/3[$.
The transformation $\TH$ translates $B_1$ onto $\TH B_1\egdef [1/3,2/3[$, and translates 
$\TH B_1$ onto $\TH^2 B_1\egdef [2/3,1[$. At this first step, $\TH$ is not yet defined on $\TH^2 B_1$. 
In general, after the $n$-th step of the construction, $X_H$ has been split into $3^n$
intervals of same length: $B_n, \TH B_n, \ldots, \TH^{3^n-1}B_n$. These intervals form a so-called 
\emph{Rokhlin tower} with base $B_n$ and height $3^n$. Such a tower is usually represented by 
putting the intervals one on top the other, the transformation $\TH$ mapping each point to the one 
exactly above. At this step, the transformation is not yet defined on $\TH^{3^n-1}B_n$. Step $n+1$ starts
by chopping the base $B_n$ into three subintervals of the same length, the first of 
which is denoted by $B_{n+1}$. The $n$-th Rokhlin tower is thus split into three 
columns, which are stacked together to get the $n+1$st tower. This amounts to mapping $\TH^{3^n-1}B_{n+1}$
onto the second piece of $B_n$ by a translation, and $\TH^{2\times 3^n-1}B_{n+1}$ onto the third piece
of $B_n$. $\TH$ is now defined everywhere except on $\TH^{3^{n+1}-1}B_{n+1}$.

The iteration of this construction for all $n\ge 1$ defines $\TH$ everywhere on $X_H$. The 
transformation obtained in this way preserves Lebesgue measure, and it is well known that
the dynamical system is ergodic. 

\subsection{The extension}

\def\pfois{\mathop{\mathord{\ldotp}\mathord{\times}}}

In order to construct the extension of $\TH$, we will now define a Rokhlin cocycle
$x\longmapsto S_x$ from $X_H$ into the group of automorphisms of $(Y,\B,\nu)$, where
$Y\egdef \{-1,1\}^{\NN}$, $\B$ is the Borel $\sigma$-algebra of $Y$, and $\nu$ is the probability measure on $Y$
which makes the coordinates independent and identically distributed, with  
$\nu(y_k=1)=\nu(y_k=-1)=1/2$ for each $k\ge 0$. 

If $y=(y_k)_{k\in\NN}\in Y$ and $0\le i\le j$, we denote by $y|_i^j$ the finite word $y_i y_{i+1}\cdots y_j$.
For each $n\ge 0$, we call \emph{$n$-block}  a word of length $2^n$ on the alphabet $\{-1,1\}$.
Le \emph{first $n$-block of $y$} is thus $y|_0^{2^n-1}$. 
If $w_1=y_0\ldots y_{2^n-1}$ and $w_2=z_0\ldots z_{2^n-1}$ are two $n$-blocks, we denote by $w_1\,w_2$ 
the $(n+1)$-block obtained by the concat\'enation of $w_1$ and $w_2$, and $w_1\pfois w_2$ the $n$-block 
defined by the termwise product of $w_1$ and $w_2$:
$$ w_1w_2\ \egdef\ y_0\ldots y_{2^n-1}z_0\ldots z_{2^n-1},\quad\mbox{and} \quad
   w_1\pfois w_2\ \egdef\ (y_0\times z_0)\ldots(y_{2^n-1}\times z_{2^n-1}). $$

For each $n\ge1$, we now define a transformation $\tn$ of Y which will be useful
for the construction of the Rokhlin cocycle. This transformation only affects
the first $n$-block of $y$ : if this first $n$-block is
$w_1w_2$ (where $w_1$ and $w_2$ are $(n-1)$-blocks), then the first $n$-block of
$\tn y$ is $w_2\,(w_1\pfois w_2)$. Coordinates with indices at least
$2^n$ of $\tn y$ remain unchanged. The two following properties 
of $\tn$ are easy to verify:
\begin{itemize}
\item $\tn$ preserves the probability $\nu$,
\item $\tn^3=\Id_Y$.
\end{itemize}

For every $x\in X_H$, we denote by $n(x)$ the smaller integer $n\ge1$ such that
$x$ does not belong to the top of tower $n$. In other words, $n(x)$ is the integer
$n\ge1$ such that $\TH x$ is at the step $n$ of the construction of $\TH$. We
then set 
$$ S_x\ \egdef\ \tau_{n(x)}\circ \tau_{n(x)-1} \circ\cdots\circ \tau_1. $$
From the properties of $\tn$, it is easy to derive that $S_x$ is always an automorphism
of $(Y,\B,\nu)$. From now on, we denote by $T$ the skew product on 
$X_H\times Y$ equipped with the product measure $\mu_H\otimes \nu$ defined by 
$$ \mathop{T}(x,y)\ \egdef \ (\TH x, \mathop{S_x}y), $$
Let $\H$ be the factor of $T$ given by the $\sigma$-algebra $\A_H\otimes \{Y,\emptyset\}$.

\subsection{Relative twofold mixing which is not threefold}

Let $n\ge 1$, and $(x,y)\in X_H\times Y$ with $x$ in the base $B_n$ of the $n$-th tower.
For each $k\ge 0$, we denote by $y^{(k)}$ the point of $Y$ defined by 
$T^k(x,y)=(\TH^k x, y^{(k)})$. From the construction of the Rokhlin cocycle, while $\TH^k x$ 
has not reached the top of tower $n$, $y$ is only transformed by some 
$\tau_j$ with $j\le n$. Therefore, in the sequence
$y^{(0)},y^{(1)},\ldots,y^{(3^n-1)}$ (corresponding to the climb of $x$ upward tower $n$), 
only the first $n$-block is modified and these modifications do not depend on 
the coordinates of $y$ with indices at least $2^n$.

We are particularly interested in the sequence $y_0^{(0)}y_0^{(1)}\ldots y_0^{(3^n-1)}$ of coordinates
with null index, which we see as a random colouring of the climb of $x$ upward tower $n$. 
From the preceding remark, this colouring only depends on the first $n$-block of $y$. 
Therefore there exists some map $\gamma_n\colon \{-1,1\}^{2^n}\to \{-1,1\}^{3^n}$
such that
$$ y_0^{(0)}y_0^{(1)}\ldots y_0^{(3^n-1)}\ =\ \gamma_n\left(y|_0^{2^n-1}\right). $$

\begin{lemma}
\label{taun+1}
  Assume further that $x$ lies in the base of the first or second column in tower $n$
  (\ie $x\in B_{n+1}$ or $x\in \TH^{3^n}B_{n+1}$). Then
  $$ y^{(3^n)}\ =\ \tnp y. $$
\end{lemma}

\begin{proof}
It is easily checked by induction on $n$, using the fact that $\tn^3=\Id_Y$.
\end{proof}

Lemma~\ref{taun+1} gives a relation between $\gamma_n$ and $\gamma_{n+1}$. Indeed, if $x$ lies 
in $B_{n+1}$, the climbing of $x$ upward tower $(n+1)$ can be seen as three successive
climbings of $x$ upward tower $n$, whose colourings are given by  $y^{(0)}=y$, $y^{(3^n)}=\tnp y$ 
and $y^{(2\times 3^n)}=\mathop{\tau_{n+1}}^2y$. It follows that the colouring of the first climbing
of $x$ upward tower $n$ is coded by the first $n$-block $y|_0^{2^n-1}$ of $y$, the colouring
of the second climbing of $x$ upward tower $n$ is coded by the second $n$-block 
$y|_{2^n}^{2^{n+1}-1}$, and the colouring of the third climbing of $x$ upward tower $n$ is 
coded by their termwise product $y|_0^{2^n-1}\pfois y|_{2^n}^{2^{n+1}-1}$. Hence, if $w$ is an
$(n+1)$-block which is the concatenation of the two $n$-blocks $w_1\,w_2$, we have
\begin{equation}
  \label{eq:recurrence}
  \gamma_{n+1}(w)\ =\ \gamma_n(w_1)\, \gamma_n(w_2)\, \gamma_n(w_1\pfois w_2).
\end{equation}
Therefore, the sequence $(\gamma_n)_{n\ge 1}$ of coding maps is entirely determined by  
$$ \gamma_1\colon a\,b\longmapsto a\,b\,(a\times b) $$
and the recurrence relation
(\ref{eq:recurrence}). The proof of the following lemma follows easily:

\begin{lemma}
\label{codage_produit}
  Let $w_1$ and $w_2$ be two $n$-blocks. Then 
  $$ \gamma_n(w_1\pfois w_2)\ =\ \gamma_n(w_1)\pfois\gamma_n(w_2). $$
\end{lemma}

From the preceding observations, we can deduce some properties of the conditional
law of the colouring process knowing $x$.

\begin{prop}
\label{colourings}
  Let $x\in X_H$ and $n\ge 1$. Let $j\ge0$ be the smallest integer such that $\TH^{-j}x\in B_{n+1}$.
  We denote by
  $C_1^n$, $C_2^n$ and $C_3^n$ the respective random colouring of the three successive climbings
  of $x$ upward tower $n$. The conditional law of $(C_1^n,C_2^n,C_3^n)$
  knowing $\H$ satisfies the following properties : 
  \begin{itemize}
  \item $C_1^n$, $C_2^n$ and $C_3^n$ are identically distributed;
  \item $C_1^n$, $C_2^n$ and $C_3^n$ are pairwise independent;
  \item $C_3^n = C_1^n \pfois C_2^n$.
  \end{itemize}
\end{prop}
\begin{proof}
Since $\H$ is a $T$-invariant $\sigma$-algebra, we can always assume to simplify the notations that $j=0$ 
(\ie $x\in B_{n+1}$). It follows from what has been seen before that
$C_1^n$, $C_2^n$ and $C_3^n$ are given respectively by $\gamma_n(y|_0^{2^n-1})$, 
$\gamma_n(y|_{2^n}^{2^{n+1}-1})$ and $\gamma_n(y|_0^{2^n-1}\pfois y|_{2^n}^{2^{n+1}-1})$.
But these three $n$-blocks $y|_0^{2^n-1}$, $y|_{2^n}^{2^{n+1}-1}$ and 
$y|_0^{2^n-1}\pfois y|_{2^n}^{2^{n+1}-1}$ are identically distributed and pairwise independent. 
Therefore, the three colourings are themselves identically distributed and pairwise independent.
The equality $C_3^n = C_1^n \pfois C_2^n$ is a straightforward consequence of Lemma~\ref{codage_produit}.
\end{proof}

It follows easily from Proposition~\ref{colourings} that the property~(\ref{eq:defmr}) characterizing
twofold mixing relatively to the factor $\H$ is true when $A$ and $B$ are measurable with
respect to a finite number of coordinates of the colouring process  $(y_0\circ T^k)_{k\in \ZZ}$. 
Indeed, in such a case we can find an integer $n$ (depending on $x$) such that $A$ and $B$ 
are measurable with respect to one of the blocks  $C_i^n$ 
($i=1,2$ or $3$) defined in the previous proposition. Then, as soon as $k\ge 3^n$, $A$ and $T^{-k}B$
are given by two blocks  $C_j^m$ (for some $m\ge n$) which are independent under the conditional law
knowing $\H$. 

Then, (\ref{eq:defmr}) extends by density to every sets $A$ and $B$ measurable with respect to 
the $\sigma$-algebra generated by $\H$ and the colouring process $(y_0\circ T^k)_{k\in \ZZ}$. 
But this $\sigma$-algebra is easily shown to be the whole $\A_H\otimes\B$, since knowing $x$ and 
$(y_0\circ T^k)_{k\in \ZZ}$ we can always recover each coordinate $y_n$, $n\in\NN$. (Details are left 
to the reader.) It follows that the system is $\H$-relatively twofold mixing.

However, the system is not $\H$-relatively threefold mixing: If $A$, $B$ and $C$ are defined by 
$$ A\ =\ B\ =\ C\ \egdef\ \{(x,y):\ y_0=1\}, $$
we have
$$ \mu(A|\H)\ =\ \mu(B|\H)\ =\ \mu(C |\H)\ =\ 1/2, $$
but for each $n\ge 1$ and each $x$ in the first column of tower $n$, 
$$ \mu(A\cap T^{-3^n}B\cap T^{-2\times3^n}C|\H)\ =\ 1/4. $$

\section{Comments and questions}

\subsection*{Joinings}

The question of the existence of a system which is twofold but not threefold mixing 
is strongly connected with the following question: Does there exist a joining of three
copies of some weakly mixing, zero-entropy dynamical system which is pairwise independent
but which is not the product measure? In \cite{leman9}, Lema\'nczyk, Mentzen and Nakada
answer positively to the \emph{relative} version of this problem: They construct 
a relatively weakly-mixing extension $T$ of an ergodic rotation $T_H$, and a 3-joining
$\lambda$ of $T$ identifying the three copies of $T_H$, which is pairwise but not 
threewise independent relative to $T_H$.
However their construction does not seem to come from an extension which is twofold
but not threefold relatively mixing.

\subsubsection*{Mixing in the base?}

The example which we have presented above can easily be modified in order to make
the dynamical system in the base weakly mixing. Indeed, we can replace the triadic Von Neumann-Kakutani 
by Chacon's transformation, whose construction is similar with the only following difference:
In each step of the construction we add a supplementary \emph{spacer} interval between
second and third column. The sequence $(h_n)$ of the heights of the successive towers thus satisfies
$h_{n+1}=3 h_n+1$. It is well known that Chacon's transformation is weakly, but not strongly, mixing. 
Defining $S_x$ in a similar way when $x$ does not lie in some spacer, and  $S_x\egdef\Id$ in
any spacer, we get the same conclusion concerning twofold but not threefold relative mixing. The lack
of threefold relative mixing is checked by considering, for $x$ in the first column of
tower $n$, $\mu(A\cap T^{-h_n}B\cap T^{-(2 h_n+1)}C|\H)$.

Then it is natural to look for a similar result with the dynamical system 
in the base strongly mixing. Indeed, it is easily shown that if  $T_H$ is mixing and 
if $T$ is $\H$-relatively mixing, then $T$ is mixing. This would give some hope to get a 
transformation that is twofold but not threefold mixing. However, there seem to be serious obstacles
to achieving the same kind of construction with a mixing base.

\subsubsection*{On the definition of relative mixing}

In the present work we have used the definition of relative mixing defined by the 
convergence to zero in probability (or equivalently in $L^1$) of the sequence
\begin{equation}
\label{suite}
\mu\left(A\cap T^{-k}B|\H\right)-\mu(A|\H)\mu(T^{-k}B|\H).
\end{equation}
An other possible definition of relative mixing is used by Rahe in his work on factors of 
Markov processes \cite{rahe1}: In this paper, a process $(x_k)_{k\in\ZZ}$ (with $x_k=x_0\circ T^k$) is said
\emph{$\H$-relatively mixing} if,
for all $A$ and $B$ measurable with respect of a finite number of cooordinates of the process $(x_k)$, 
the convergence of~(\ref{suite}) to zero holds almost surely.

The difference between these two definitions is discussed in a recent work of Rudolph \cite{rudol9},
where it is shown that there exists a system which is relatively mixing with respect to one of
his factors in the $L^1$ sense, but not in the almost-sure sense. Rudolph also shows that checking
almost-sure convergence of~(\ref{suite}) to zero for a dense class of subsets $A$ and $B$ (as in Rahe's
definition) implies that the same convergence holds for \emph{every} $A$ and $B$.

It is not difficult to see that, for the example we presented here, the same
results concerning twofold and threefold relative mixing hold if we replace $L^1$ convergence
by almost-sure convergence. 

\providecommand{\bysame}{\leavevmode\hbox to3em{\hrulefill}\thinspace}

\end{document}